# Periodic solutions of a class of second-order non-autonomous differential equations

Jia Ruan

*College of Mathematics and Statistics, Chongqing University, Chongqing 401331, China*

**Abstract**: This paper investigates the dynamical behavior of periodic solutions for a class of second-order non-autonomous differential equations. First, based on the Lyapunov-Schmidt reduction method for finite-dimensional functions, the corresponding bifurcation function is constructed, and it is proven that the system possesses at least one *T*-periodic solution. Second, a two-timing method is employed to perform perturbation analysis on the original equation. By separating the fast and slow time scales, an explicit expression for the approximate *T*-periodic solution is derived. Furthermore, for the stability of the system under parametric excitation, the bifurcation characteristics near the first instability tongue are revealed through perturbation expansion and eigenvalue analysis. Additionally, the Ince-Strutt stability diagram is plotted to illustrate the stability boundaries.

**Keywords:** Periodic solution; Perturbed differential equation; Lyapunov-Schmidt reduction method; Two-timing method; Ince-Strutt stability diagram

## 1 Introduction

In recent years, with advancements in modeling across fields such as engineering vibrations, circuit design, population evolution, and automatic control, there has been an increasing demand for research on the solutions and properties of nonlinear differential equations. In particular, the periodic solutions of second-order nonlinear differential equations have attracted significant attention. First, the existence and stability of periodic solutions constitute fundamental theoretical issues that require priority investigation, as they form the basis for further discussions on the properties of solutions. Second, during the modeling process, small perturbations that are not precisely quantified may exhibit transient or persistent characteristics. By introducing such disturbance terms into perturbed differential equations, the realism of the model can be enhanced, and the properties of the original solutions may be substantially altered. Therefore, in-depth analysis of the stability, periodicity, and asymptotic behavior of solutions to perturbed differential equations is critical for understanding the system's sensitivity dependence on minor perturbations.

The classical Mathieu equation is a special case of Hill's equation [1], introduced by the French mathematician Mathieu when analyzing the vibration modes of a stretched membrane with elliptic clamped boundaries. Many mechanical and electrical models can be expressed by the Mathieu equation, among which the most popular mechanical model is the parametrically excited pendulum, also known as the Kapitza pendulum. Figure 1.1 illustrates this concept:

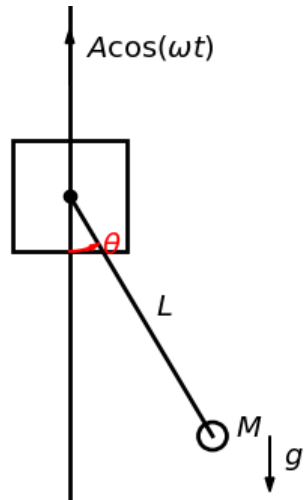

Fig.1 Kapitza pendulum

This system consists of a massless rigid rod and a small, heavy bob. Its pivot point experiences periodic motion in the vertical direction. The governing differential equation for this system [2] is

$$\theta'' + \left(\frac{g}{L} - \frac{A\omega^2}{L}\cos(\omega t)\right)\sin\theta = 0, \tag{1}$$

where $M, g, L, A, \omega$ represent, respectively, the mass of the bob, gravitational acceleration, rod length, amplitude of the vertical motion, and excitation frequency. When $|\theta|$ is small, equation (1) can be approximated by the Taylor expansion as

$$\theta'' + \left(\frac{g}{L} - \frac{A\omega^2}{L}\cos(\omega t)\right)\theta = 0, \tag{2}$$

which represents the linear Kapitza pendulum.

A forced inverted pendulum [3] can be generally described by

$$x'' + [\alpha + \beta q(t)]x = f(t), \tag{3}$$

where $q(t+T) = q(t)$ is the parametrically excited function and $f(t)$ is the forced $kT$-periodic function with $k \in Z$, both of them have zero mean value.

Hence, for the nonlinear Kapitza pendulum with external forcing, one may write

$$\theta'' + \left(\frac{g}{L} - \frac{A\omega^2}{L}\cos\omega t\right)\theta = f(t),$$

where $f(t)$ is a $kT$-periodic function and $T = 2\pi/\omega$.

In the fields of science, physics, and engineering, a specialized vibration equation featuring a cubic nonlinear term has been derived. This equation integrates the cubic nonlinear characteristics inherent to the classical Mathieu equation and the Duffing equation, and is thus designated as the Mathieu-Duffing equation.

Mohsen Azimi [2] studied a class of undamped Mathieu-Duffing equations with cubic nonlinearities, expressed as

$$x'' + (\omega_n^2 + F_p \cos(\omega_p t))x + \alpha x^3 = 0, \tag{4}$$

Where $\omega_n, \omega_p, F_p, \alpha$ are constants representing natural frequency, parametric frequency, parametric forcing amplitude and nonlinear stiffness. Perturbing the time-varing and nonlinear terms in equation (4) results in the following equation.

$$x'' + (\omega_n^2 + \varepsilon \cos(\omega_p t))x + \varepsilon \alpha x^3 = 0. \tag{5}$$

Applying the method of strained parameter, Azimi analyzed stability and bifurcation near the first and the second instability tongues for this equation.

In this paper, we introduce a

$$x'' + (\omega_n^2 + \varepsilon \cos(\omega_p t))x + \varepsilon \alpha x^3 = \varepsilon f(t), \tag{6}$$

where $f(t)$ satisfies Dirichlet conditions, has periodic $T = 2\pi/\omega_p$, and has zero mean over one periodic.

Our primary contributions are :

1.**Bifurcation of *T*-periodic solutions under Resonance condition.**

When $\omega_n = \omega_p = \omega$, the equation (6) reduces to

$$x'' + (\omega^2 + \varepsilon \cos(\omega t))x + \varepsilon \alpha x^3 = \varepsilon f(t). \tag{7}$$

By applying the Lyapunov-Schmidt reduction method, we establish the corresponding bifurcation function and prove the existence of *T*-periodic solutions.

**2. Two-Timing Method and Stability Analysis.**

Using the two-timing method, we derive an approximate *T*-periodic solution for equation (7) by separating fast and slow time scales. We further analyze the stability and bifurcation behavior of the equation (8) in the vicinity of the first instability tongue, exploring how the interplay of parametric and external forcing shapes the system's dynamic responses.

# 2 The existence of periodic solutions

In this chapter, we imploy the Lyapunov-Schmidt reduction method [4] to study the bifurcation of *T*-periodic solutions for a class of second-order non-autonomous equation.

## 2.1 Basic lemma

We consider the problem of the bifurcation of *T*-periodic solutions from the differention system

$$x'(t) = F_0(t,x) + \varepsilon F_1(t,x) + \varepsilon F_2(t,x,\varepsilon), \tag{8}$$

with $\varepsilon = 0$ to $\varepsilon \neq 0$ sufficienty small. The functions $F_0, F_1 : \mathbb{R} \times \Omega \to \mathbb{R}^n$ and $F_2 : \mathbb{R} \times \Omega \times (-\varepsilon_f, \varepsilon_f) \to \mathbb{R}^n$ are $C^2$ in their variables and *T*-periodic in the first variable, and $\Omega$ is an open subset of $\mathbb{R}^n$. One of the main assumptions is that the unperturbed system

$$x'(t) = F_0(t,x) \tag{9}$$

has an open subset of $\Omega$ fulfilled of periodic solutions.

Let $x = x(t,z,\varepsilon)$ be the solution of system (8) such that $x(0,z,\varepsilon) = z$. We write the linearization of the unperturbed system (9) along a periodic solution $x(t,z,0)$ as

$$y' = D_x F_0(t, x(t,z,0))y. \tag{10}$$

Let $Y(\cdot, z)$ be some fundamental matrix solution of (10). Then we have the following lemmas.

**Lemma 1**. [5] We assume that there exists an open set $V$ with $\overline{V} \subset \Omega$ and such that for each $z \in \overline{V}$, $x(t,z,0)$ is *T*-periodic and we consider the function $f_1 : \overline{V} \to \mathbb{R}^n$ given by

$$f_1(z) = \int_0^T Y^{-1}(t,z) F_1(t, x(t,z,0)) dt.$$

If there exists $a \in V$ with $f_1(a) = 0$ and $\det((df_1/d\alpha)(a)) \neq 0$, then there exists a $T$-periodic solution $\varphi(\cdot, \varepsilon)$ of system (9) such that $\varphi(0, \varepsilon) \to a$ as $\varepsilon \to 0$.

## 2.2 Existence of a *T*-Periodic solution

This section employs the Lyapunov-Schmidt reduction method to analyze the bifurcation problem of *T*-periodic solutions for equation (7). For clarity, we first transform the equation into a first-order differential system. Let $x' = y$, then equation (7) becomes:

$$\begin{cases} x' = y, \\ y' = -\omega^2 x + \varepsilon(f(t) - \alpha x^3 - \cos(\omega t)x), \end{cases} \tag{11}$$

Where $\omega > 0$, $\alpha \neq 0$ and $\varepsilon$ is a small parameter. $f(t)$ is a *T*-periodic function where $T = 2\pi/\omega$.

Before proceeding, we verify whether the perturbed system (11) satisfies the prerequisites for the Lyapunov-Schmidt reduction. Consider the unperturbed system:

$$\begin{cases} x' = y, \\ y' = -\omega^2 x. \end{cases} \tag{12}$$

The general solution of system (12) with initial conditions $x(0) = x_0$, $y(0) = y_0$ is

$$\begin{cases} x(t) = \dfrac{\omega x_0 \cos(\omega t) + y_0 \sin(\omega t)}{\omega}, \\ y(t) = -\omega x_0 \sin(\omega t) + y_0 \cos(\omega t). \end{cases}$$

For and initial conditions $(x_0, y_0) \neq (0,0)$, all solutions of system (12) are $2\pi/\omega$-periodic. Thus, we can invoke Lemma 1 from section 2.1 to handle the perturbed differential system (11). Based on this, we establish the existence of *T*-periodic solutions for a class of second-order non-autonomous differential equations.

**Theorem 1.** Suppose $f(t)$ in (11) is a *T*-periodic function, $T = 2\pi/\omega$, satisfying the Dirichlet conditions and having zero average (i.e., $\int_0^T f(t)dt = 0$). Then for sufficiently small $\varepsilon$, there exists an initial condition $x(0) = x_0^*$, $y(0) = y_0^*$, such that system (11) admits at least one *T*-periodic solution $(x(t,\varepsilon), y(t,\varepsilon))$, where

$$(x_0^*, y_0^*) = \left( -\sqrt[3]{\frac{4a_1^3}{3\alpha\left(a_1^2 + b_1^2\right)}}, -\sqrt[3]{\frac{4b_1^3 \omega^3}{3\alpha\left(a_1^2 + b_1^2\right)}} \right).$$

Moreover, as $\varepsilon \to 0$, the solution

$$(x(t,\varepsilon), y(t,\varepsilon)) \to \left( \frac{\omega x_0^* \cos(\omega t) + y_0^* \sin(\omega t)}{\omega}, -\omega x_0^* \sin(\omega t) + y_0^* \cos(\omega t) \right).$$

**Proof:** Rewrite system (11) in the form

$$u'(t) = F_0(t,u) + \varepsilon F_1(t,u),$$

where

$$u = \begin{pmatrix} x \\ y \end{pmatrix}, \quad F_0(t,u) = \begin{pmatrix} y \\ -\omega^2 x \end{pmatrix},$$

$$F_1(t,u) = \begin{pmatrix} 0 \\ f(t) - \alpha x^3 - \cos(\omega t)x \end{pmatrix}. \tag{13}$$

We denote the linearization of (12) by

$$\frac{d}{dt}\begin{pmatrix} y_1 \\ y_2 \end{pmatrix} = \begin{pmatrix} 0 & 1 \\ -\omega^2 & 0 \end{pmatrix}\begin{pmatrix} y_1 \\ y_2 \end{pmatrix},$$

which has two linearly independent solutions

$$Y_{11} = \begin{pmatrix} \cos(\omega t) \\ -\omega \sin(\omega t) \end{pmatrix}, \quad Y_{12} = \frac{1}{\omega}\begin{pmatrix} \sin(\omega t) \\ \omega \cos(\omega t) \end{pmatrix}.$$

Therefore, the fundamental matrix solution is

$$Y(t,z) = \begin{pmatrix} \cos(\omega t) & \dfrac{\sin(\omega t)}{\omega} \\ -\omega \sin(\omega t) & \cos(\omega t) \end{pmatrix}$$

with its inverse

$$Y^{-1}(t,z) = \begin{pmatrix} \cos(\omega t) & -\dfrac{\sin(\omega t)}{\omega} \\ \omega \sin(\omega t) & \cos(\omega t) \end{pmatrix}. \tag{14}$$

To establish the existence of $T$-periodic solutions, we compute the bifurcation via Lemma 1.

$$f_1(z) = \int_0^T Y^{-1}(t,z)F_1(t,u(t,z,0))dt, \tag{15}$$

where

$$z = \begin{pmatrix} x(0) \\ y(0) \end{pmatrix} = \begin{pmatrix} x_0 \\ y_0 \end{pmatrix},$$

$$u(t,z,0) = \begin{pmatrix} \dfrac{\omega x_0 \cos(\omega t) + y_0 \sin(\omega t)}{\omega} \\ -\omega x_0 \sin(\omega t) + y_0 \cos(\omega t) \end{pmatrix}. \tag{16}$$

Substituting (13), (14) and (16) into (15), we obtain

$$f_1(x_0, y_0) = \int_0^{\frac{2\pi}{\omega}} \begin{pmatrix} \cos(\omega t) & -\dfrac{\sin(\omega t)}{\omega} \\ \omega \sin(\omega t) & \cos(\omega t) \end{pmatrix}\begin{pmatrix} 0 \\ f(t) - \alpha x^3 - \cos(\omega t)x \end{pmatrix}\bigg|_{x = \frac{\omega x_0 \cos(\omega t) + y_0 \sin(\omega t)}{\omega}} dt. \tag{17}$$

The calculation is devided into three parts, let

$$M_1(x_0, y_0) = \int_0^{\frac{2\pi}{\omega}} \begin{pmatrix} \cos(\omega t) & -\frac{\sin(\omega t)}{\omega} \\ \omega \sin(\omega t) & \cos(\omega t) \end{pmatrix} \begin{pmatrix} 0 \\ f(t) \end{pmatrix} \bigg|_{x = \frac{\omega x_0 \cos(\omega t) + y_0 \sin(\omega t)}{\omega}} dt, \quad (18)$$

$$M_2(x_0, y_0) = \int_0^{\frac{2\pi}{\omega}} \begin{pmatrix} \cos(\omega t) & -\frac{\sin(\omega t)}{\omega} \\ \omega \sin(\omega t) & \cos(\omega t) \end{pmatrix} \begin{pmatrix} 0 \\ \alpha x^3 \end{pmatrix} \bigg|_{x = \frac{\omega x_0 \cos(\omega t) + y_0 \sin(\omega t)}{\omega}} dt, \quad (19)$$

$$M_3(x_0, y_0) = \int_0^{\frac{2\pi}{\omega}} \begin{pmatrix} \cos(\omega t) & -\frac{\sin(\omega t)}{\omega} \\ \omega \sin(\omega t) & \cos(\omega t) \end{pmatrix} \begin{pmatrix} 0 \\ \cos(\omega t) x \end{pmatrix} \bigg|_{x = \frac{\omega x_0 \cos(\omega t) + y_0 \sin(\omega t)}{\omega}} dt. \quad (20)$$

Therefore, equation (17) can be expressed as

$$f_1(x_0, y_0) = M_1(x_0, y_0) - M_2(x_0, y_0) - M_3(x_0, y_0). \quad (21)$$

Since $f(t)$ satisfying the Dirichlet conditions and has zero mean over one periodic, it can be expanded in a convergent Fourier series [6]

$$f(t) = \sum_{n=1}^{\infty} (a_n \cos(n\omega t) + b_n \sin(n\omega t)), \quad (22)$$

where $a_1$ and $b_1$ are not all equal to 0.

Substituting (22) into (18), we obtain

$$M_1(x_0, y_0) = \int_0^{\frac{2\pi}{\omega}} \begin{pmatrix} -\frac{\sin(\omega t)}{\omega} \left( \sum_{n=1}^{\infty} (a_n \cos(n\omega t) + b_n \sin(n\omega t)) \right) \\ \cos(\omega t) \left( \sum_{n=1}^{\infty} (a_n \cos(n\omega t) + b_n \sin(n\omega t)) \right) \end{pmatrix} dt.$$

$$= \int_0^{\frac{2\pi}{\omega}} \begin{pmatrix} -\frac{1}{\omega} \left( a_1 \sin(\omega t) \cos(\omega t) + b_1 \sin^2(\omega t) \right) \\ a_1 \cos^2(\omega t) + b_1 \sin(\omega t) \cos(\omega t) \end{pmatrix} dt$$

$$= \frac{\pi}{\omega^2} \begin{pmatrix} -b_1 \\ \omega a_1 \end{pmatrix}. \quad (23)$$

When $n$ is even, we have

$$\int_0^{\frac{\pi}{2}} \sin^n x \, dx = \int_0^{\frac{\pi}{2}} \cos^n x \, dx = \frac{(n-1)!!}{n!!} \cdot \frac{\pi}{2}.$$

And when $k$ and $l$ are positive odd numbers,

$$\int_0^{\frac{2\pi}{\omega}} \sin^k(\omega t) \cos^l(\omega t) dt = 0.$$

Therefore, equation (19) can be expresses as

$$M_2(x_0, y_0) = \frac{\alpha}{\omega^3} \int_0^{\frac{2\pi}{\omega}} \begin{pmatrix} -\frac{1}{\omega}\sin(\omega t)(\omega x_0 \cos(\omega t) + y_0 \sin(\omega t))^3 \\ \cos(\omega t)(\omega x_0 \cos(\omega t) + y_0 \sin(\omega t))^3 \end{pmatrix} dt = \begin{pmatrix} -\dfrac{3\pi\alpha y_0 (\omega^2 x_0^2 + y_0^2)}{4\omega^5} \\ \dfrac{3\pi\alpha x_0 (\omega^2 x_0^2 + y_0^2)}{4\omega^3} \end{pmatrix}. \quad (24)$$

Through a similar analysis, equation (20) can be reduced to

$$M_3(x_0, y_0) = 0. \quad (25)$$

Substituting (23), (24) and (25) into (21), we obtain

$$f_1(x_0, y_0) = \begin{pmatrix} -\dfrac{\pi}{\omega^2} b_1 + \dfrac{3\pi\alpha y_0 (\omega^2 x_0^2 + y_0^2)}{4\omega^5} \\ \dfrac{\pi}{\omega} a_1 - \dfrac{3\pi\alpha x_0 (\omega^2 x_0^2 + y_0^2)}{4\omega^3} \end{pmatrix}.$$

After obtaining the bifurcation function of system (11), we proceed to find the zeros of this bifurcation function. Let

$$f_{11}(x_0, y_0) = -\frac{\pi}{\omega^2} b_1 + \frac{3\pi\alpha y_0 (\omega^2 x_0^2 + y_0^2)}{4\omega^5},$$

$$f_{21}(x_0, y_0) = \frac{\pi}{\omega} a_1 - \frac{3\pi\alpha x_0 (\omega^2 x_0^2 + y_0^2)}{4\omega^3}.$$

The zeros of the bifurcation function must satisfy the following system of equations

$$\begin{cases} f_{11}(x_0, y_0) = 0, \\ f_{21}(x_0, y_0) = 0. \end{cases}$$

Solving these equations yields

$$\begin{cases} x_0^* = \sqrt[3]{\dfrac{4a_1^3}{3\alpha(a_1^2 + b_1^2)}}, \\ y_0^* = \sqrt[3]{\dfrac{4b_1^3 \omega^3}{3\alpha(a_1^2 + b_1^2)}}. \end{cases}$$

Additionally, the Jacobian determinant is computed as

$$\det\left(\frac{\partial(f_{11}, f_{21})}{\partial(x_0, y_0)}\right) = \begin{vmatrix} \dfrac{3\pi\alpha x_0 y_0}{2\omega^3} & \dfrac{3\pi\alpha(\omega^2 x_0^2 + 3y_0^2)}{4\omega^5} \\ -\dfrac{3\pi\alpha(y_0^2 + 3\omega^2 x_0^2)}{4\omega^3} & -\dfrac{3\pi\alpha x_0 y_0}{2\omega^3} \end{vmatrix}$$

$$= -\frac{9\pi^2 \alpha^2 x_0^2 y_0^2}{4\omega^6} + \frac{9\pi^2 \alpha^2 [(\omega^2 x_0^2 + 3y_0^2)(y_0^2 + 3\omega^2 x_0^2)]}{16\omega^8}$$

$$= \frac{27\pi^2 \alpha^2}{16\omega^8}(y_0^2 + \omega^2 x_0^2)^2.$$

Since $\alpha \neq 0$, $\omega > 0$, it is evident that

$$\det\left(\frac{\partial(f_{11}, f_{21})}{\partial(x_0, y_0)}\bigg|_{(x_0, y_0)=(x_0^*, y_0^*)}\right) \neq 0$$

By Lemma 1, the system (11) with initial conditions $x(0) = x_0^*$, $y(0) = y_0^*$ admits a $T$-periodic solution $(x(t,\varepsilon), y(t,\varepsilon))$, where $(x_0^*, y_0^*)$ are given in (3.27). Furthermore, this periodic solution satisfies

$$(x(t,\varepsilon), y(t,\varepsilon)) \to \left(\frac{\omega x_0^* \cos(\omega t) + y_0^* \sin(\omega t)}{\omega}, -\omega x_0^* \sin(\omega t) + y_0^* \cos(\omega t)\right)$$

as $\varepsilon \to 0$. □

## 3 Existence of approximate periodic solutions

Perturbation methods are widely used in mathematics and physics to address problems that are difficult to solve directly. In this chapter, we employ the two-timing method [7], which offers better approximation, to study the existence of approximate periodic solutions for a class of second-order non-autonomous differential equations.

**Theorem 2.** Consider the nonlinear equation with a small parameter

$$x'' + (\omega^2 + \varepsilon \cos(\omega t))x + \varepsilon \alpha x^3 = \varepsilon \sum_{n=1}^{\infty}(a_n \cos(n\omega t) + b_n \sin(n\omega t)), \qquad (26)$$

where $\omega > 0$, $\alpha \neq 0$, $|\varepsilon| \ll 1$, and the fourier coefficients $a_1$ and $b_1$ of the driving term are not all zero. Then, Equation (26) admits an approximate $T$-periodic solution with $T = 2\pi/\omega$. The zeroth-order approximate solution is

$$x(t) = \frac{\omega x_0^* \cos(\omega t) + y_0^* \sin(\omega t)}{\omega} + O(\varepsilon),$$

where

$$(x_0^*, y_0^*) = \left(\sqrt[3]{\frac{4a_1^3}{3\alpha(a_1^2 + b_1^2)}}, \sqrt[3]{\frac{4b_1^3 \omega^3}{3\alpha(a_1^2 + b_1^2)}}\right).$$

**Proof:** Introduce two time scales: the fast time scale $T_0 = t$ and the slow time scale $T_1 = \varepsilon t$. Using the chain rule, we derive

$$\frac{\partial x}{\partial t} = \frac{\partial x}{\partial T_0} + \varepsilon \frac{\partial x}{\partial T_1},$$

and

$$\frac{\partial^2 x}{\partial t^2} = \frac{\partial^2 x}{\partial T_0^2} + 2\varepsilon \frac{\partial^2 x}{\partial T_0 \partial T_1} + \varepsilon^2 \frac{\partial^2 x}{\partial T_1^2}. \qquad (27)$$

Substituting (27) into (26) yields

$$\frac{\partial^2 x}{\partial T_0^2} + 2\varepsilon \frac{\partial^2 x}{\partial T_0 \partial T_1} + \varepsilon^2 \frac{\partial^2 x}{\partial T_1^2} + (\omega^2 + \varepsilon \cos(\omega t))x + \varepsilon \alpha x^3 = \varepsilon \sum_{n=1}^{\infty}(a_n \cos(n\omega t) + b_n \sin(n\omega t)). \quad (28)$$

Expand the solution of (26) in terms of the small parameter

$$x(t) = x_0(T_0, T_1) + \varepsilon x_1(T_0, T_1) + \varepsilon^2 x_2(T_0, T_1) + \cdots. \quad (29)$$

Substitute (29) into (28) and neglect higher-order terms $O(\varepsilon^2)$. Separating orders of $\varepsilon$, the zeroth-order equation is

$$\frac{\partial^2 x_0}{\partial T_0^2} + \omega^2 x_0 = 0. \quad (30)$$

The first-order equation is

$$\frac{\partial^2 x_1}{\partial T_0^2} + \omega^2 x_1 = -2\frac{\partial^2 x_0}{\partial T_0 \partial T_1} + \sum_{n=1}^{\infty}(a_n \cos(n\omega T_0) + b_n \sin(n\omega T_0)) - \alpha x_0^3 - \cos(\omega T_0) x_0. \quad (31)$$

The solution to (30) corresponds to simple harmonic motion

$$x_0(T_0, T_1) = A(T_1)\cos(\omega T_0 + \psi(T_1)), \quad (32)$$

Substituting (32) into (31), we analyze the resonance terms. By eliminating secular terms (terms leading to unbounded solutions), we derive the amplitude and phase modulation equations

$$\begin{cases} 2\omega \dfrac{dA}{dT_1} + a_1 \sin\psi(T_1) + b_1 \cos\psi(T_1) = 0, \\ 2\omega A(T_1)\dfrac{d\psi}{dT_1} + a_1 \cos\psi(T_1) - b_1 \sin\psi(T_1) - \dfrac{3\alpha A^3(T_1)}{4} = 0. \end{cases} \quad (33)$$

From Equation (32), it can be observed that the equilibrium point $(A_0, \psi_0)$ of system (33) corresponds to a periodic motion of system (26). Below, we seek the equilibrium point $(A_0, \psi_0)$ of system (33) by setting

$$\begin{cases} \dfrac{dA}{dT_1} = 0, \\ \dfrac{d\psi}{DT_1} = 0. \end{cases}$$

From equation (33), we obtain

$$\begin{cases} a_1 \sin\psi_0 + b_1 \cos\psi_0 = 0, \\ a_1 \cos\psi_0 - b_1 \sin\psi_0 = \dfrac{3\alpha A_0^3}{4}. \end{cases} \quad (34)$$

By introducing the variable substitution

$$x_0^* = A_0 \cos\psi_0, \quad y_0^* = -A_0 \sin\psi_0. \quad (35)$$

Equation (34) yields

$$\begin{cases} b_1 x_0^* = a_1 y_0^*, \\ a_1 x_0^* + b_1 y_0^* = \dfrac{3\alpha}{4}\left((x_0^*)^2 + (y_0^*)^2\right)^2. \end{cases} \quad (36)$$

Is is noteworthy that equation (34) implies

$$a_1^2 + b_1^2 = \frac{9\alpha^2 A_0^6}{16}$$

Since $a_1$ and $b_1$ are not both zero, $A_0 \neq 0$. Therefore, $x_0^*$ and $y_0^*$ are not both zero. Solving equation (36), we obtain

$$\begin{cases} x_0^* = \sqrt[3]{\dfrac{4 a_1^3}{3\alpha\left(a_1^2 + b_1^2\right)}}, \\ y_0^* = \sqrt[3]{\dfrac{4 b_1^3}{3\alpha\left(a_1^2 + b_1^2\right)}}. \end{cases} \quad (37)$$

Substituting equations (37) into (32), the zeroth-order solution becomes

$$\begin{aligned} x_0(T_0) &= A_0 \cos(\omega T_0 + \psi_0) \\ &= A_0 \cos(\omega T_0)\cos\psi_0 - A_0 \sin(\omega T_0)\sin\psi_0 \\ &= x_0^* \cos(\omega T_0) + y_0^* \sin(\omega T_0) \\ &= \sqrt[3]{\dfrac{4 a_1^3}{3\alpha\left(a_1^2 + b_1^2\right)}} \cos(\omega T_0) + \sqrt[3]{\dfrac{4 b_1^3}{3\alpha\left(a_1^2 + b_1^2\right)}} \sin(\omega T_0). \end{aligned} \quad (38)$$

Thus, based on equations (33) and (38), the solution to equation (26) can be expressed as

$$x(t) = \frac{\omega x_0^* \cos(\omega t) + y_0^* \sin(\omega t)}{\omega} + O(\varepsilon),$$

where

$$(x_0^*, y_0^*) = \left( \sqrt[3]{\dfrac{4 a_1^3}{3\alpha\left(a_1^2 + b_1^2\right)}},\ \sqrt[3]{\dfrac{4 b_1^3}{3\alpha\left(a_1^2 + b_1^2\right)}} \right).$$

Hence, we have derived an approximate $T$-periodic solution for equation (26). □

It is worth noting that as $\varepsilon \to 0$, the solution of equation (26) satisfies

$$x(t) \to \frac{\omega x_0^* \cos(\omega t) + y_0^* \sin(\omega t)}{\omega}.$$

## 4. Stability and bifurcation analysis

Inspired by Mohsen Azimi's work on parametric resonance and nonlinear dynamics

[2], this section examines the stability of a nonlinear differential equation with periodic forcing:

$$\frac{\partial^2 x}{\partial t^2} + \left(\omega_n^2 + \varepsilon \cos(\omega_p t)\right)x + \varepsilon \alpha x^3 = \varepsilon \sum_{n=1}^{\infty}(a_n \cos(n\omega_p T_0) + b_n \sin(n\omega_p T_0)), \quad (39)$$

where $\omega_n$ is the natural frequency, $\omega_p$ is the parametric frequency, $\alpha \neq 0$, $\varepsilon$ is a small parameter, and $a_1, b_1$ are not all zero. The stability of solutions depends critically on parameter combinations and is visualized via Ince-Strutt diagrams [8].

Using the two-timing method, we expand the solution as

$$\frac{\partial^2 x}{\partial t^2} = \frac{\partial^2 x}{\partial T_0^2} + 2\varepsilon \frac{\partial^2 x}{\partial T_0 \partial T_1} + \varepsilon^2 \frac{\partial^2 x}{\partial T_1^2}. \quad (40)$$

Substituting equation (40) into (39), we obtain

$$\frac{\partial^2 x}{\partial T_0^2} + 2\varepsilon \frac{\partial^2 x}{\partial T_0 \partial T_1} + \varepsilon^2 \frac{\partial^2 x}{\partial T_1^2} + (\omega_n^2 + \varepsilon \cos(\omega_p T_0))x + \varepsilon \alpha x^3 = \varepsilon \sum_{n=1}^{\infty}(a_n \cos(n\omega_p T_0) + b_n \sin(n\omega_p T_0)). \quad (41)$$

The solution to equation (39) is expanded in terms of the small parameter as

$$x(t) = x_0(T_0, T_1) + \varepsilon x_1(T_0, T_1) + O(\varepsilon^2). \quad (42)$$

Using the perturbation method [9], the natural frequency in equation (39) is expanded as

$$\omega_n^2 = \frac{1}{4}\omega_p^2 + \varepsilon \omega_1 + O(\varepsilon^2). \quad (43)$$

Substituting equations (42) and (43) into equation (39) and neglecting higher-order terms $O(\varepsilon^2)$, we separate the orders of $\varepsilon$. The zeroth-order equation is

$$\frac{\partial^2 x_0}{\partial T_0^2} + \frac{1}{4}\omega_p^2 x_0 = 0. \quad (44)$$

The first-order equation is

$$\frac{\partial^2 x_1}{\partial T_0^2} + \frac{1}{4}\omega_p^2 x_1 = -2\frac{\partial^2 x_0}{\partial T_0 \partial T_1} - x_0 \cos(\omega_p T_0) - \omega_1 x_0 - \alpha x_0^3$$

$$+ \sum_{n=1}^{\infty}\left(a_n \cos(n\omega_p T_0) + b_n \sin(n\omega_p T_0)\right). \quad (45)$$

The system corresponding to equation (44) represents simple harmonic motion, with the general solution

$$x_0(T_0, T_1) = A(T_1)\cos\left(\frac{1}{2}\omega_p T_0 + \psi(T_1)\right), \quad (46)$$

Substituting equation (46) into equation (45), we analyze the resonance terms.

By eliminating secular terms (terms leading to unbounded solutions), we derive the amplitude and phase modulation equations

$$\begin{cases} \omega_p \dfrac{dA}{dT_1} - \dfrac{1}{2} A(T_1)\sin(2\psi(T_1)) = 0, \\ \omega_p A(T_1) \dfrac{d\psi}{dT_1} - \omega_1 A(T_1) - \dfrac{3}{4}\alpha A^3(T_1) - \dfrac{1}{2} A(T_1)\cos(2\psi(T_1)) = 0. \end{cases} \quad (47)$$

System (47) has a trivial solution $A = 0$. For non-trivial solutions $A \neq 0$, Equation (47) simplifies to:

$$\begin{cases} \dfrac{dA}{dT_1} = \dfrac{A(T_1)\sin(2\psi(T_1))}{2\omega_p}, \\ \dfrac{d\psi}{dT_1} = \dfrac{\omega_1 + \dfrac{3}{4}\alpha A^2(T_1) + \dfrac{1}{2}\cos(2\psi(T_1))}{\omega_p}. \end{cases} \quad (48)$$

It is noteworthy that the equilibrium points $(A^*, \psi^*)$ of system (48) correspond to periodic motions of system (39).

To find the equilibrium points $(A^*, \psi^*)$, we set

$$\begin{cases} \dfrac{dA}{dT_1} = 0, \\ \dfrac{d\psi}{dT_1} = 0. \end{cases}$$

Solving these yields

$$\begin{cases} \psi = 0, \dfrac{\pi}{2}, \pi, \dfrac{3}{2}\pi, \\ A^2 = -\dfrac{4}{3\alpha}\left[\dfrac{\cos(2\psi)}{2} + \omega_1\right]. \end{cases} \quad (49)$$

For $\alpha < 0$ (the case $\alpha > 0$ is analogous and omitted here), the condition $A^2 > 0$ implies

(1) If $\psi = 0$ or $\pi$, $\omega_1$ must satisfy $\omega_1 > -\dfrac{1}{2}$.

(2) If $\psi = \dfrac{\pi}{2}$ or $\dfrac{3\pi}{2}$, $\omega_1$ must satisfy $\omega_1 > \dfrac{1}{2}$.

Thus, the transition curve are

$$\omega_n^2 = \dfrac{1}{4}\omega_p^2 \pm \dfrac{1}{2}\varepsilon + O(\varepsilon^2).$$

To analyze the stability of the equilibrium points $(A^*, \psi^*)$, we introduce polar coordinates

$$M(T_1) = A(T_1)\cos\psi(T_1), N(T_1) = -A(T_1)\sin\psi(T_1). \qquad (50)$$

Using equation (50), equation (46) becomes

$$x_0(T_0,T_1) = A(T_1)\cos\left(\frac{1}{2}\omega_p T_0 + \psi(T_1)\right) = M(T_1)\cos\left(\frac{1}{2}\omega_p T_0\right) + N(T_1)\sin\left(\frac{1}{2}\omega_p T_0\right).$$

In polar coordinates,

$$A^2(T_1) = M^2(T_1) + N^2(T_1), \tan\psi(T_1) = -\frac{N(T_1)}{M(T_1)}. \qquad (51)$$

From equation (51), we derive

$$AA' = MM' + NN', \psi' = \frac{M'N - N'M}{M^2 + N^2}. \qquad (52)$$

Substituting equation (52) into equation (48), we obtain

$$\begin{cases} M' = \dfrac{1}{\omega_p}\left[\omega_1 N - \dfrac{1}{2}N + \dfrac{3\alpha N}{4}(M^2 + N^2)\right], \\ N' = \dfrac{1}{\omega_p}\left[-\omega_1 M - \dfrac{1}{2}M - \dfrac{3\alpha M}{4}(M^2 + N^2)\right]. \end{cases} \qquad (53)$$

The Jacobian matrix of system (53) at point $(M, N)$ is

$$J = \begin{pmatrix} \dfrac{3\alpha MN}{2\omega_p} & \dfrac{\omega_1 - \dfrac{1}{2} + \dfrac{3}{4}\alpha(M^2 + 3N^2)}{\omega_p} \\ \dfrac{-\omega_1 - \dfrac{1}{2} - \dfrac{3}{4}\alpha(3M^2 + N^2)}{\omega_p} & -\dfrac{3\alpha MN}{2\omega_p} \end{pmatrix}. \qquad (54)$$

From Equations (49) and (50), the equilibrium points of system (53) are as follows

1. For $\omega_1 < -\dfrac{1}{2}$, there is one equilibrium point $(M_1, N_1)$.
2. For $-\dfrac{1}{2} < \omega_1 < \dfrac{1}{2}$, there are three equilibrium points $(M_1, N_1), (M_2, N_2), (M_3, N_3)$.
3. For $\omega_1 > \dfrac{1}{2}$, there are five equilibrium points $(M_1, N_1), (M_2, N_2), (M_3, N_3)$, $(M_4, N_4), (M_5, N_5)$,

where

$$(M_1, N_1) = (0,0), (M_2, N_2) = \left(\sqrt{-\frac{4}{3\alpha}\left(\omega_1 + \frac{1}{2}\right)}, 0\right), (M_3, N_3) = \left(-\sqrt{-\frac{4}{3\alpha}\left(\omega_1 + \frac{1}{2}\right)}, 0\right),$$

$$(M_4, N_4) = \left(0, \sqrt{-\frac{4}{3\alpha}\left(\omega_1 - \frac{1}{2}\right)}\right), (M_5, N_5) = \left(0, -\sqrt{-\frac{4}{3\alpha}\left(\omega_1 - \frac{1}{2}\right)}\right).$$

To analyze the stability of these equilibrium points, we employ the following lemma

**Lemma 4.1** [2]. The local stability of an equilibrium point of system (53) is determined by the eigenvalues of its Jacobian matrix (54). Since the trace of $J$ is zero $\mathrm{tr}(J)=0$, the eigenvalues of $J$ are a pair of opposites. The stability criteria are:

1. If $\det(J)>0$, the eigenvalues are purely imaginary, and the equilibrium is a stable center.

2. If $\det(J)<0$, the eigenvalues are real with opposite signs, and the equilibrium is an unstable saddle.

The stability of the five equilibrium points is analyzed as follows:

For $(M_1, N_1)$, $\det(J) = \dfrac{4\omega_1^2 - 1}{4\omega_p}$. By Lemma 2, $(M_1, N_1)$ is unstable when $|\omega_1| < \dfrac{1}{2}$ and stable when $|\omega_1| > \dfrac{1}{2}$.

For $(M_2, N_2)$ and $(M_3, N_3)$, $\det(J) = \dfrac{2\omega_1 + 1}{\omega_p^2}$. These points exist for $\omega_1 > -\dfrac{1}{2}$, where $\det(J) > 0$, indicating stability.

For $(M_4, N_4)$ and $(M_5, N_5)$, $\det(J) = \dfrac{-2\omega_1 + 1}{\omega_p^2}$. These points exist for $\omega_1 > \dfrac{1}{2}$, where $\det(J) < 0$, indicating instability.

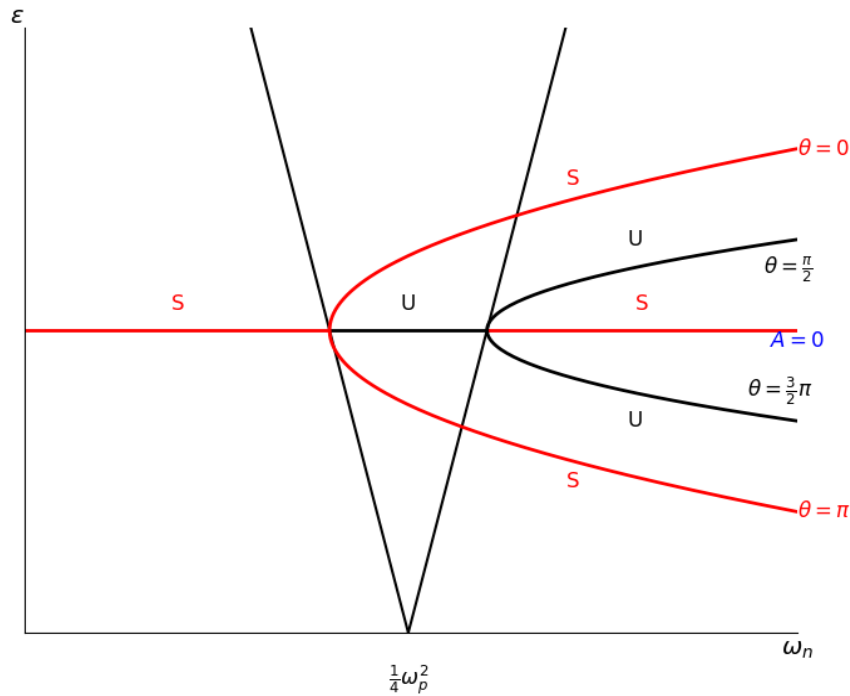

Fig.2 Strutt-Ince diagram

Conclusion: As shown in Figure 2, when the nonlinear stiffness $a<0$ in equation (39), increasing $\omega_n$ while keeping $\varepsilon$ constant leads to the following bifurcations:

1. A Supercritical pitchfork bifurcation occurs upon crossing the left transition curve. The trivial solution $(M_1, N_1) = (0,0)$ becomes unstable, and two new stable equilibria (centers) $(M_2, N_2)$ 和 $(M_3, N_3)$ emerge.

2. A Subcritical pitchfork bifurcation occurs upon crossing the right transition curve. The trivial solution regains stability, and two new unstable saddle points $(M_4, N_4)$ 和 $(M_5, N_5)$ emerge.

Remark: At the transition curves $\det(J) = 0$, the stability of equilibria typically changes [10], as illustrated in Figure 2.